\documentclass{amsart}
\usepackage{amsmath,amscd,amsthm,amsfonts,a4wide,mathrsfs}
\usepackage[all]{xy}

\theoremstyle{plain}
\newtheorem{theorem}                 {Theorem}      [section]
\newtheorem{proposition}  [theorem]  {Proposition}
\newtheorem{corollary}    [theorem]  {Corollary}
\newtheorem{lemma}        [theorem]  {Lemma}

\theoremstyle{definition}
\newtheorem{example}      [theorem]  {Example}
\newtheorem{remark}       [theorem]  {Remark}

\numberwithin{equation}{section}

\def \g{\mathfrak{g}}
\def \k{\mathfrak{k}}
\def \p{\mathfrak{p}}

\def \SU#1{{\bf SU}(#1)}
\def \S{{\bf S}}
\def \U#1{{\bf U}(#1)}
\def \su#1{\mathfrak{su}(#1)}
\def \s{\mathfrak{s}}
\def \u#1{\mathfrak{u}(#1)}
\def \SLC#1{{\bf SL}_{#1}(\cn)}
\def \slc#1{\mathfrak{sl}_{#1}(\cn)}

\def \L{\mathscr{L}} 
\def \d{\mathrm{d}}
\def \dstar{\delta}
\def \ep{\langle\cdot,\cdot\rangle}
\def \epn{(\cdot,\cdot)}
\def \ip#1#2{\langle#1,#2\rangle}
\def \abs#1{\lvert\,#1\,\rvert}
\def \span{\mathrm{span}}

\DeclareMathOperator{\volume}{\ast 1}
\DeclareMathOperator{\trace}{trace}
\DeclareMathOperator{\Hom}{Hom}
\DeclareMathOperator{\Ad}{Ad}

\DeclareMathOperator{\musicd}{\sharp}
\DeclareMathOperator{\musicu}{\flat}

\def \rn{\mathbb{R}}
\def \cn{\mathbb{C}}

\newcommand{\ee}{\mathscr{E}} 
\newcommand{\R}{{\mathbb{R}}} 
\newcommand{\Z}{{\mathbb{Z}}} 
 
\newcommand{\C}{{\mathbb{C}}}

\newcommand{\CP}{{\mathbb{C}}{{P}}}

\newcommand{\beq}{\begin{equation}} 
\newcommand{\eeq}{\end{equation}} 
\newcommand{\bea}{\begin{eqnarray}} 
\newcommand{\eea}{\end{eqnarray}} 
\newcommand{\ben}{\begin{eqnarray*}} 
\newcommand{\een}{\end{eqnarray*}} 
\newcommand{\ra}{\rightarrow}

\newcommand{\thet}{{\vartheta}}

\newcommand{\spec}{{\rm spec\, }}

\newcommand{\ol}{\overline}

\newcommand{\id}{{\rm Id}} 
\newcommand{\etavec}{\mbox{\boldmath${\eta}$}}

\allowdisplaybreaks

\begin{document}\larger[2]\setlength{\baselineskip}{1.0\baselineskip}

\title[On the Strong Coupling Limit of the Faddeev-Hopf Model]{On the Strong Coupling Limit of the Faddeev-Hopf Model}
\author[J.M.~Speight]{J.M.~Speight}
\thanks{The second author was supported by the Swedish 
Research Council (623-2004-2262)}
\subjclass[2000]{58E99, 81T99}
\address{School of Mathematics, University of Leeds, Leeds, LS2 9JT} 
\email{speight@maths.leeds.ac.uk }
\author[J.M.L.~Svensson]{M.~Svensson}
\address{School of Mathematics, University of Leeds, Leeds, LS2 9JT} 
\email{M.Svensson@leeds.ac.uk }

\begin{abstract} 
The variational calculus for the Faddeev-Hopf model on a general 
Riemannian domain, with general K\"ahler target space, is
studied in the strong coupling
limit. In this limit, the model has key similarities with pure
Yang-Mills theory, namely conformal invariance in dimension $4$
and an infinite dimensional symmetry group.
The first and second variation formulae are calculated and several
examples of stable solutions are obtained. In particular, it is proved that
all immersive solutions are stable. Topological lower energy
bounds are found in dimensions $2$ and $4$.
An explicit description of the spectral behaviour 
of the Hopf map $S^3\ra S^2$ is given, 
and a conjecture of Ward 
concerning the stability of this map in the full Faddeev-Hopf model 
is proved.
\end{abstract}

\maketitle

\section{Introduction}

Theoretical physics has long been a rich source of geometrically interesting 
and natural variational problems. The Yang-Mills equations, of deep
significance for the differential topology of 4-manifolds \cite{donkro}, and
the Yang-Mills-Higgs equations, which have led to interesting results in
hyperk\"ahler geometry \cite{atihit}, both originated in elementary 
particle physics. Harmonic map theory, while not originating in theoretical
physics, has found many applications in high energy and condensed matter 
physics, with physicists frequently contributing genuinely new insights. 

The purpose of this paper is to present a systematic study of a variational
problem arising in the so-called Faddeev-Hopf (or
Faddeev-Skyrme) model \cite{faddeev}, originally proposed as a 
model of quark confinement (among other phenomena) in
high energy physics. Let $M$ be some Riemannian manifold, representing 
physical space, and $N$ a K\"ahler manifold, the
target space, with K\"ahler form $\omega$. The model has a single field
$\phi:M\ra N$, the energy functional (or action functional, in the case
where $M$ is euclideanized spacetime after Wick rotation) being
$$
\ee(\phi)=\frac12\int_M(|\d\phi|^2+\alpha|\phi^*\omega|^2),
$$
$\alpha\geq 0$ being a coupling constant. The model of original interest has
$M=\R^3$, $N=S^2$. The weak coupling limit of this model, 
$\alpha=0$ has of course been intensively 
studied: it is the harmonic map problem. This is conformally invariant if
$M$ has dimension $2$. By contrast, we shall study the strong coupling limit,
$\alpha\ra\infty$, or more precisely, the variational problem for the
energy functional
$$
E(\phi)=\lim_{\alpha\ra\infty}\alpha^{-1}\ee(\phi)=\frac12\int_M
|\phi^*\omega|^2.
$$
This does not seem to have received systematic study in either the 
theoretical
physics or differential geometry communities. It has been studied in the
specific case $M=\R\times S^3$ (with a Lorentzian metric, actually)
and $N=S^2,\C$ or the hyperbolic plane by de Carli and Ferreira 
\cite{decfer}.
It has some important similarities with pure Yang-Mills 
theory. It is invariant under an infinite dimensional group of symmetries,
the group of symplectic diffeomorphisms of $N$, rather as Yang-Mills
theory is invariant under gauge transformations. It is also,
as we will demonstrate, conformally
invariant if $M$ has dimension $4$. Both these facts were known to
de Carli and Ferreira in the specific context they studied. The most 
interesting situation physically is when $M=S^4$, interpreted as the
conformal compactification of $\R^4$. Nontrivial solutions in this case may
receive the physical interpretation of instantons in the strong coupling
limit of the Faddeev-Hopf model in $(3+1)$ dimensions, just as critical 
points
of the Yang-Mills functional on $S^4$ are interpreted as pure gauge-theory
instantons. Such solutions have profound effects on the quantized version
of the field theory \cite[ch10]{raj}. 

Our motivation for studying this variational problem is twofold. First, 
simple curiosity prompts us to ask what the geometric character of the 
variational calculus for this functional is. We will see that both the
first and second variation formulae can be given elegant and natural
geometric formulations from which strong results quickly follow. For example,
we will show that all immersive solutions are stable, and that there are
no non-vacuum (i.e., $E>0$)
immersive solutions in the case $M=S^4$, for any choice of target space. 
Second, we
hope that studying one term in the Faddeev-Hopf model in isolation will
give valuable insight into the finite coupling model. Indeed, we will
identify a large class of critical points of $E$ which are also harmonic
maps, and hence critical points of $\ee$ for all $\alpha$. In particular,
we are able to prove a stability conjecture
of Ward concerning the full Faddeev-Hopf model on $M=S^3$ \cite{war}. The 
functional $E$ also arises as one term in the so-called baby Skyrme models
studied by Zakrzewski and collaborators \cite{zak}, and our results should 
find applications in these models too. 

The rest of the paper is structured as follows. In section \ref{firvar}
we carefully define the functional $E$, prove that it is conformally 
invariant
in dimension $4$, derive the first variation formula (Euler-Lagrange 
equation)
for $\phi$ and construct some interesting explicit solutions. In section
\ref{sub} we consider submersive solutions in particular, identifying
a large class of critical submersions which are also harmonic. In section
\ref{stable} we obtain topological lower bounds on $E$ when $M$ has dimension
$2$ or $4$, and derive the second variation formula in the general case. 
The results are used to prove the stability of several interesting solutions.
Finally, in section \ref{hopf} 
the variational calculus for the projection $G\ra G/K$ onto
a Hermitian symmetric space is developed in general, and the results used
to show that the Hopf map $S^3\ra S^2$ in particular is stable. A proof of
Ward's conjecture quickly follows from this.

\section{The First Variation}\label{firvar}

In this section we assume that $(M^m,g)$ is a compact, oriented
Riemannian manifold of dimension $m$. For any vector bundle $E$ over $M$, 
we denote by $\Gamma(E)$ the space of sections of $E$. 

Recall that the metric $g$ on $M$ induces a (pointwise) metric on the 
bundle of $p$-forms on $M$, defined by 
$$
\ip{\alpha}{\beta}=\frac{1}{p!}\sum_{i_1,\dots,i_p=1}^m\alpha(e_{i_1},
\dots,e_{i_p})\beta(e_{i_1},\dots,e_{i_p}),
$$ 
where $e_1,\dots,e_m$ is a local orthonormal frame on $M$. By using the 
Hodge $\ast$-operator, we get the relation 
$$
\alpha\wedge\ast\beta=\ip{\alpha}{\beta}\volume
$$ 
for any two $p$-forms $\alpha$ and $\beta$; here $\volume$ is the volume 
element on $M$. Integrating this inner product over $M$ gives a global 
$L^2$-product 
$$
\ip{\alpha}{\beta}_{L^2}=\int_M\ip{\alpha}{\beta}\volume=
\int_M\alpha\wedge\ast\beta\qquad(\alpha,\beta\in\Gamma(\wedge^pT^*M),
$$ 
with corresponding norm $\|\alpha\|^2_{L^2}=\ip{\alpha}{\alpha}_{L^2}$. 
With respect to this $L^2$-product, the exterior differentiation operator 
$$
\d:\Gamma(\wedge^pT^*M)\to\Gamma(\wedge^{p+1}T^*M)
$$ 
has the adjoint operator 
$$
\dstar:\Gamma(\wedge^{p}T^*M)\to\Gamma(\wedge^{p-1}T^*M),\quad \dstar\alpha
=(-1)^{m+mp+1}\ast\d\ast\alpha.
$$ 

We will also be using the \emph{musical isomorphisms} on $M$ which are 
defined as follows: 
$$
\musicu:\Gamma(TM)\to\Gamma(T^*M),\quad \musicu X=g(X,\cdot),\quad 
\musicd=\musicu^{-1}.
$$ 

Let $(N^n,h,J)$ be a compact K\"ahler manifold of real dimension $n$ and 
with K\"ahler form $\omega=h(J\cdot,\cdot)$. For a smooth map 
$\phi:M\to N$ we define the \emph{energy
functional} 
$$
E(\phi)=\frac12\|\phi^*\omega\|^2_{L^2}=\frac12\int_M\phi^*\omega\wedge\ast\phi^*\omega.
$$ 
Any map $\phi$ for which $E(\phi)=0$, the minimum possible,
 will be called a \emph{vacuum solution}
or \emph{vacuum} of the theory. Clearly $\phi$ is a vacuum if and only if
$\phi^*\omega=0$ everywhere, that is, if $\phi$ is isotropic.

We begin our investigation of $E(\phi)$ by verifying that it is, like the 
Yang-Mills functional, invariant under conformal changes of $g$ if $M$ has 
dimension $4$. 
\begin{proposition}\label{conf}
 Assume that $\phi:M\to N$ is a map from a $4$-dimensional Riemannian 
manifold to a K\"ahler manifold. Then the functional $E(\phi)$, and 
therefore also the Euler-Lagrange equation for $\phi$, is invariant under 
conformal changes of the metric on $M$. 
\end{proposition}

\begin{proof} Assume that $g$ is a Riemannian metric on $M$ and that 
$\tilde g=\lambda^2 g$, where $\lambda$ is some positive function on $M$. 
Denote by $\ep_{\tilde g}$ the metric induced by $\tilde g$ on 2-forms and 
let $\tilde\ast 1$ be the corresponding volume element. Then 
$$
\tilde\ast 1=\lambda^4\volume\text{ and }\ip{\alpha}{\alpha}_{\tilde g}
=\lambda^{-4}\ip{\alpha}{\alpha},
$$ 
for any 2-form $\alpha$ on $M$. Hence the form $\ip{\alpha}{\alpha}\volume$ 
remains unchanged.  
\end{proof}

Next we derive the Euler-Lagrange equation for $E(\phi)$ in the general case.
\begin{proposition}\label{proposition:first} For a smooth variation 
$\phi_t:M\to N$ of $\phi$ with variational vector field 
$X\in\Gamma(\phi^{-1}TN)$, we have 
$$
\frac{d}{dt}E(\phi_t)\big|_{t=0}=
\int_M\omega(X,\d\phi(\musicd\dstar\phi^*\omega))\volume.
$$ 
\end{proposition}

For the proof, let us recall the following simple result.
\begin{lemma}[Homotopy Lemma] Let $M$ and $N$ be two manifolds and 
$$
\phi_t:M\to N
$$ 
a smooth family of maps. For any closed $2$-form $\eta$ on $N$ we have 
$$
\frac{\partial}{\partial t}\phi_t^*\eta=
\d\big(\phi_t^*\iota(\frac{\partial\phi_t}{\partial t})\eta\big).
$$
\end{lemma} 

Here $\iota$ denotes the interior product. For a proof of this lemma see, 
e.g., \cite[p49]{Ee-Le:83}.

\begin{proof}[Proof of Proposition \ref{proposition:first}] 
It is obvious that $\dfrac{\partial}{\partial t}$ commutes with the Hodge 
$\ast$-operator on $\wedge^* T^*M$. By the Homotopy Lemma we have 
\begin{equation*}
\begin{split}
\frac12\frac{\partial}{\partial t}\big|_{t=0}\phi^*_t\omega\wedge\ast\phi^*_t\omega=&\frac12(\frac{\partial}{\partial t}\phi^*_t\omega)\wedge\ast\phi^*_t\omega+\frac12\phi^*_t\omega\wedge\ast(\frac{\partial}{\partial t}\phi^*_t\omega)\\
=&\frac12\d\phi^*\iota_X\omega\wedge\ast\phi^*_t\omega+\frac12\phi^*\omega\wedge\ast\phi^*\iota_X\omega\\
=&\ip{\d\phi^*\iota_X\omega}{\phi^*\omega}\volume.
\end{split}
\end{equation*}
Therefore 
$$
\frac{\partial}{\partial t}E(\phi_t)\big|_{t=0}=
\int_M\ip{\phi^*\iota_X\omega}{\dstar\phi^*\omega}\volume,
$$ 
and 
$$
\ip{\phi^*\iota_X\omega}{\dstar\phi^*\omega}=
\phi^*\iota_X\omega(\musicd\dstar\phi^*\omega)=
\omega(X,\d\phi(\musicd\dstar\phi^*\omega)).
$$
This proves the proposition. 
\end{proof}

\begin{corollary}\label{cor:ker} The map $\phi:M\to N$ is a critical point 
for the functional if and only if 
$$
\musicd\dstar\phi^*\omega\in\ker\d\phi
$$ 
everywhere on $M$. 
\end{corollary}

\begin{example} Assume that $M=N$ and $\phi:N\to N$ is the identity map. 
Then $\phi^*\omega=\omega$, and $\dstar\omega=0$. Hence $\phi$ is a 
critical point for the functional. 
\end{example}

\begin{remark}\label{harm} Assume that $\phi$ is smooth and is
immersive on a dense set, that is, for all $x$ in a dense subset of $M$,
the differential 
$$\d\phi_x:T_xM\to T_{\phi(x)}N$$ 
is injective. By the 
Corollary, if $\phi$ is a critical point of the functional, then 
the 1-form $\dstar\phi^*\omega$ vanishes almost everywhere, and hence 
vanishes everywhere by continuity. Hence
the 2-form $\phi^*\omega$ is co-closed. Since it is obviously closed, we 
see that \emph{an almost everywhere immersive map is a critical point of 
the functional if and only if $\phi^*\omega$ is a harmonic $2$-form}. 

In particular, when $H^2(M,\rn)=0$, the only immersive critical points 
defined on $M$ are \emph{isotropic} immersions, i.e., maps for which 
$\phi^*\omega=0$. As previously remarked, such a map has $E(\phi)=0$, and 
hence is a 
vacuum solution of the field theory. The set
of vacuum solutions of this theory is unusually rich.
\end{remark}

\begin{example}
The map 
$
\phi:S^4\to\cn P^4,
$ 
defined as the 2-fold covering by $S^4$ of $\rn P^4$ followed by the natural 
embedding of $\rn P^4$ to $\cn P^4$ is clearly isotropic, hence a vacuum.
An interesting question is whether this vacuum is path connected, through
vacua, to the trivial vacuum $\phi=$ constant. One suspects not.

Certainly the set of vacua may fail to be path connected. 
Consider the zero section of
$TS^n$, equipped with the Stenzel metric \cite{stenzel}. This is manifestly
an isotropic immersion $i:S^n\ra TS^n$, hence a vacuum. Given any
smooth map $\phi:S^n\ra S^n$, the composition $i\circ\phi$ is still 
isotropic, hence a vacuum. But if the degree of $\phi$ is not unity, then
$i$ and $i\circ\phi$ are not even homotopic, much less path connected through
isotropic maps. Hence the set of vacua in the case $M=S^n$ and $N=TS^n$
is not path connected.
 
\end{example}

Of primary physical interest, given their physical interpretation as 
instantons, are smooth anisotropic critical points on $M=S^4$
which minimize $E(\phi)$ within their homotopy 
class. 
In particular, one would like a smooth anisotropic
minimizer in the nontrivial
class of $\pi_4(S^2)$, a \emph{pure Faddeev-Hopf instanton}.  
In fact, it remains an open question whether smooth anisotropic
critical points exist on $S^4$ at all, for any choice of target space. 
A standard starting point for finding special solutions is the use
of symmetry reduction.
A fundamental difficulty in exploiting symmetry reduction
is raised by Remark
\ref{harm}:  symmetry reductions of the variational problem on $S^4$
tend to produce only maps which are either trivial or immersive. The best we
have managed is a smooth solution mapping the punctured hemisphere into
$\CP^2$, two
copies of which can be glued together to give a continuous map $S^4
\ra \CP^2$ which, away from
the poles and the equator of $S^4$,
is smooth and satisfies the field equation, and has finite total energy. 
This map
is constructed in the next example.

\begin{example} 
Let $M=S^4$ and $N=\CP^2$. The twice punctured 4-sphere
is conformally equivalent to the cylinder $\R\times \SU{2}$ given the
metric 
$$
g=dt^2+\sigma_1^2+\sigma_2^2+\sigma_3^2,
$$
where $\sigma_i$ are the usual left-invariant one forms on $\SU{2}$. So we
may seek critical points on $\R\times\SU{2}$ satisfying appropriate
boundary conditions as $|t|\ra\infty$. Let $I=(a,b)\subset\R$ be any open 
interval and $Q$ be the Banach manifold of, for example, $C^2$ maps
$I\times\SU{2}\ra\CP^2$. Then $E:Q\ra\R$ is $C^1$, and there is a natural
action of $\SU{2}\times V_4$ on $Q$ given by
\ben
\phi(t,X)&\stackrel{(U,1)}{\ra}&
\left[\begin{array}{cc}U& \\& 1 \end{array}\right]\phi(t,UX)\\
\phi(t,X)&\stackrel{(1,P_1)}{\ra}&
\left[\begin{array}{cc}X& \\& 1 \end{array}\right]
\left[\begin{array}{ccc}1&&\\&-1& \\&& 1 \end{array}\right]
\left[\begin{array}{cc}X^{-1}& \\& 1 \end{array}\right]
\phi(t,UX)\\
\phi(t,X)&\stackrel{(1,P_2)}{\ra}&
\ol{\phi(t,\ol{X})},
\een
where $P_1,P_2$ generate the Viergruppe $V_4=\{1,P_1,P_2,P_1P_2\}$.

The set $Q_0$ of fixed points of this action consists of maps of the
form
$$
\phi(t,X)=\left[\begin{array}{cc}X&\\&1\end{array}\right][\alpha(t):0:1],
$$
where $\alpha:I\ra S^1=\R\cup\{\infty\}$. Since $\SU{2}\times V_4$ is
compact, we may apply the principle of symmetric
criticality \cite{palais} to deduce that any critical point of $E|_{Q_0}$ 
is automatically a critical point of $E$. Routine calculation shows that
$$
E|_{Q_0}(\alpha)=\pi^2\int_I
\frac{2\alpha^2\dot{\alpha}^2}{(1+\alpha^2)^4}+
\frac{\alpha^4}{2(1+\alpha^2)^2}\ \d t,
$$
where $\dot{}$ denotes differentiation with respect to $t$. 
This may be thought of as the action of a one-dimensional
Lagrangian mechanical system. By invariance under $t$ translation, 
all solutions $\alpha(t)$ conserve the 
quantity
$$
H=\frac{\alpha^2\dot{\alpha}^2}{(1+\alpha^2)^4}-
\frac{\alpha^4}{4(1+\alpha^2)^2}.
$$
If we wish $\phi$ to extend to $S^4$ we should insist that
$\alpha,\dot{\alpha}\ra 0$ as $t\ra \infty$, so only solutions with
$H=0$
are of interest. The $H=0$ level curve in the 
$(\alpha,\dot{\alpha})$ plane is,
$
4\dot{\alpha}^2=\alpha^2(1+\alpha^2)^2$, whence one finds solutions
$\alpha_+:(0,\infty)\ra\R$ and $\alpha_-:(-\infty,0)\ra\R$,
$$
\alpha_+(t)=\frac{1}{\sqrt{e^t-1}},\qquad
\alpha_-(t)=\frac{-1}{\sqrt{e^{-t}-1}}.
$$
Gluing these together gives a continuous map $S^4\ra\CP^2$,
$$
\phi(t,X)=\left[\begin{array}{cc}X&\\&1\end{array}\right][1:0:
\frac{t}{|t|}\sqrt{e^{|t|}-1}],
$$
of total energy $E=\pi^2$,
which away from $t=\pm\infty$ and $0\times\SU{2}$,
is smooth and solves the field equation. Clearly this solution fails to
be globally smooth.
\end{example}

\begin{remark} One can find global, smooth, anisotropic, finite energy
critical maps on $\R^4$ if one equips it with a metric outside the
Euclidean conformal class. An example is given in the next
section, Example \ref{funky}.
\end{remark}

\begin{remark} De Carli and Ferreira have constructed an ingenious
symmetry reduction in the case $M=\R\times S^3$ (Lorentzian in their
original version, but the 
same reduction works in the Riemannian case) and $N=S^2$, by imposing
invariance under a $T^2\times\Z_2$ group of symmetries
\cite{decfer}. This reduces the
field equation, not to a nonlinear ODE as in the example above,
but to a \emph{linear} elliptic PDE for a single real function
on $\R^2$. Unfortunately, only solutions which remain bounded on the
whole of $\R^2$ give globally well defined maps $\phi$ on $S^4$, and no such
nontrivial solutions exist in the Riemannian case.
\end{remark}

\section{Critical Submersions}\label{sub}

In light of Corollary \ref{cor:ker} and Remark \ref{harm} it is natural to 
seek
critical maps in the case where the dimension of $M$ exceeds that of $N$.
In particular, there exists a large number of interesting critical
submersions. Such solutions automatically have non-vanishing energy, so
are not vacua. We begin with the simple example of projection on a
(possibly) warped product, then reformulate the first variation formula
in a way better suited to submersions.

\begin{example}\label{ex:warped} Assume that $(P,k)$ and $(N,h)$ are two 
compact Riemannian manifolds and that $f:P\to\rn$ is a positive, smooth 
function. The \emph{warped product} of $(P,k)$ and $(N,h)$ by $f$ is the 
manifold $P\times N$ with the Riemannian metric 
$$g=k+f^2h.$$  
Assume further that $(N,h,\omega)$ is K\"ahler. Then the projection map
$\phi:P\times N\ra N$ onto the second coordinate is critical.

To prove this, let $\ast_P$ and $\ast_N$ be the Hodge star operators on $P$ 
and $N$, respectively, so that the volume form on $P$ is $\ast_P1$. Then, 
as is easily seen, 
$$
\ast\phi^*\omega=f^{2(n-2)}\ast_P1\wedge\ast_N\omega=
\frac{f^{2(n-2)}}{(n-1)!}\ast_P1\wedge\omega^{n-1},
$$ 
where $n=\dim_\cn N$. As $f$ is a function only on $P$, 
$\d f\wedge\ast_P1=0$. Furthermore, $\ast_P1$ and $\omega^{n-1}$ are 
obviously closed. Hence
$\dstar\phi^*\omega=0$.
\end{example}

\begin{example} 
\label{funky}
One can use projection on a product to
construct global, smooth, finite energy solutions on
$M=\R^4$ if one equips it with a metric outside the Euclidean conformal 
class. For example, let 
$$
g=(1+x_3^2+x_4^2)^2(dx_1^2+dx_2^2)+(1+x_1^2+x_2^2)^2(dx_3^2+dx_4^2).
$$
Then $(\R^4,g)$ is complete, has infinite volume and is Ricci positive
with bounded scalar curvature. It is conformally equivalent to $S^2_\times
\times S^2_\times$ where $S^2_\times$ is the punctured unit sphere.
In terms of a stereographic coordinate on
$S^2_\times$ projected from the puncture, 
the equivalence is $x\equiv (x_1+ix_2,x_3+ix_4)$. 
Hence,
by Proposition \ref{conf} and Example \ref{ex:warped}, the projection map
$\R^4\ni x\mapsto x_1+ix_2\in S^2$ is critical and has energy $8\pi^2$.
\end{example}

To proceed further in our analysis of critical submersions, let us denote 
by $\nabla$ both the Levi-Civita connexion on $TM$ and on $TN$. Recall that 
the connexion on $TN$ induces a connexion $\nabla^\phi$ on $\phi^{-1}TN$. 
This connexion, together with the Levi-Civita connexion on $TM$, induces a 
connexion on $\Hom(TM,\phi^{-1}TN)$, which we also denote by $\nabla$. The 
\emph{second fundamental form} of $\phi$ is the covariant derivative of 
$\d\phi$: 
$$
\nabla\d\phi(X,Y)=\nabla^\phi_X\d\phi(Y)-
\d\phi(\nabla_XY)\qquad(X,Y\in\Gamma(TM)).
$$ 
The map $\phi$ is said to be 
\emph{totally geodesic} if its second fundamental form vanishes. The 
\emph{tension field} of $\phi$ is the trace of the second fundamental form: 
$$
\tau(\phi)=\trace\nabla\d\phi=\sum_{i=1}^m\nabla\d\phi(e_i,e_i).
$$ 
The map $\phi$ is said to be a \emph{harmonic map} if its tension field 
vanishes. 

To simplify our calculations, let us fix a point $x\in M$ and an 
orthonormal frame $\{e_i\}_{i=1}^m$ which is \emph{normal} at $x$, i.e., 
$$
\nabla_{e_i}e_j(x)=0
$$ 
for all $i,j$. Assuming that all calculations take place at the point 
$x\in M$, we can rewrite the Euler-Lagrange equations in the following way:
\begin{equation*}\label{eq:FirstVarNew}
\begin{split}
\d\phi(\musicd\dstar\phi^*\omega)=&
\sum_{j=1}^m\d\phi(\dstar\phi^*\omega(e_j)e_j)\\
=&-\d\phi\bigg(\sum_{i,j=1}^m\nabla_{e_i}(\phi^*\omega(e_i,e_j)e_j)\bigg)\\
=&-\sum_{i,j=1}^m\bigg(e_i\omega(\d\phi(e_i),\d\phi(e_j))\bigg)\d\phi(e_j)\\
=&\sum_{j=1}^m\bigg(\sum_{i=1}^m\omega(\nabla\d\phi(e_i,e_j),\d\phi(e_i))-
\omega(\tau(\phi),\d\phi(e_j))\bigg)\d\phi(e_j).
\end{split}
\end{equation*}
We thus define 
\begin{equation}\label{eq:S}
S(\phi)=\sum_{j=1}^m\bigg(\sum_{i=1}^m\omega(\nabla\d\phi(e_i,e_j),
\d\phi(e_i))-\omega(\tau(\phi),\d\phi(e_j))\bigg)\d\phi(e_j);
\end{equation} 
by our calculation, the map $\phi$ is a critical point if and only if 
$S(\phi)=0$.

\begin{example}\label{Example:Riemann} Assume that $\phi$ is a Riemannian  
submersion; thus, at each point $x\in M$, the differential  
$$
\d\phi_x:T_xM\to T_{\phi(x)}N
$$ 
maps the space $(\ker\d\phi_x)^\perp\subseteq T_xM$ isometrically onto 
$T_{\phi(x)}N$. 
We first demonstrate that  
$$
\nabla\d\phi(X,Y)=0\qquad(X,Y\in(\ker\d\phi_x)^\perp).
$$
Any vector field $X$ on $N$ can be written as $\d\phi(\hat X)$ for some 
vector field $\hat X$ on $M$ taking values in $(\ker\d\phi)^\perp$, and we 
can always find a local frame for $(\ker\d\phi)^\perp$ of vector fields of 
the form $\hat X$. Thus, it is enough to show that 
$$
\nabla\d\phi(\hat X,\hat Y)=0\qquad(X,Y\in\Gamma(TN)).
$$
Denoting by $\nabla^M$ and $\nabla^N$ the Levi-Civita connexions on $M$ and 
$N$, respectively, we have for $X,Y,Z\in\Gamma(TN)$, 
\begin{equation*}
\begin{split}
g(\nabla^M_{\hat X}\hat Y,\hat Z)=&\frac{1}{2}\big\{\hat Xg(\hat Y,\hat Z)+
\hat Yg(\hat X,\hat Z)-\hat Zg(\hat X,\hat Y)\\
&+g([\hat X,\hat Y],\hat Z)+g([\hat Z,\hat X],\hat Y)-g([\hat Y,\hat Z],
\hat X)\big\}\\
=&\frac{1}{2}\big\{Xh(Y,Z)+Yh(X,Z)-Zh(X,Y)\\
&+h([X,Y],Z)+h([Z,X],Y)-h([Y,Z],X)\big\}\\
=&h(\nabla^N_XY,Z)=g(\widehat{\nabla^N_XY},\hat Z).
\end{split}
\end{equation*}
Thus, 
$$
\nabla\d\phi(\hat X,\hat Y)=\nabla^N_XY-\d\phi(\nabla^M_{\hat X}{\hat Y})=0
\qquad(X,Y\in\Gamma(TN)).
$$ 

Locally, we can choose an orthonormal frame $e_1,\dots,e_m$ for $TM$ with 
the property that $e_1,\dots,e_{m-n}$ is a local frame for $\ker\d\phi$ and 
$e_{m-n+1},\dots,e_m$ is a local frame for $(\ker\d\phi)^\perp$. Then 
$$
S(\phi)=-\sum_{j=m-n+1}^m\omega(\tau(\phi),\d\phi(e_j))\d\phi(e_j).
$$ 
Thus, $\phi$ is critical if and only if $\phi$ is harmonic. 
In fact, using the same local frame for $TM$ gives 
$$
\tau(\phi)=-\d\phi(\sum_{j=1}^{m-n}\nabla_{e_j}e_j)=-\d\phi(H),
$$ 
where $H$ is the mean curvature vector of the fibres of $\phi$. We conclude 
that \emph{a Riemannian submersion is a critical point if and only if it 
has minimal fibres, and thus is a harmonic morphism}, see \cite{Ba-Wo:03}. 
Note that such a map, being harmonic, is automatically a critical point 
of the full Faddeev-Hopf functional for every value of the coupling, not
just the infinite coupling limit.

For example, the natural projection 
$$
\phi:S^{2n+1}\to\cn P^n,\quad \phi(z)=[z]\qquad(z\in S^{2n+1}
\subset\cn^{n+1})
$$ 
is a Riemannian submersion with minimal, even totally geodesic, fibres.
 
\end{example}

\section{The Second Variation and Stability}\label{stable}

In this section we calculate the second variation of the energy functional. 
Assume that $\phi:M\to N$ is a critical point of the functional. We define 
the \emph{Hessian of $E$ at $\phi$} as 
$$
H_\phi(X,Y)=\frac{\partial^2}{\partial s\partial t}\big|_{s=t=0}
E(\phi_{s,t});
$$ 
here $\phi_{s,t}$ is a 2-parameter variation of $\phi$ with 
$$
X=\partial_t\phi_{s,t}|_{s=t=0}\ \text{ and }\ 
Y=\partial_s\phi_{s,t}|_{s=t=0}.
$$ 
Clearly $H_{\phi}$ is a symmetric, bi-linear form on $\Gamma(\phi^{-1}TN)$. 
The map $\phi$ is said to be \emph{stable} if 
$$
H_{\phi}(X,X)\geq0\qquad(X\in\Gamma(\phi^{-1}TN));
$$ 
the \emph{index} of $\phi$ is the dimension of the largest subspace on 
which $H_{\phi}$ is negative. 

Clearly, any map which minimizes the energy within its homotopy class 
is a stable critical point. In some situations it is easy to give lower 
bounds for the energy.

\begin{proposition}\label{bound1} Let $M=M^2$ be a surface. Then
$$
E(\phi)\geq\frac{1}{2\mathrm{Vol}(M)}\big(\int_M\phi^*\omega\big)^2.
$$
Note that the right hand side is a homotopy invariant.
If $\phi$ attains this  lower bound  then
$\phi$ is either isotropic (so $E(\phi)=0$) or has no critical points.
\end{proposition}

\begin{proof}
Since $M$ is a surface, $\phi^*\omega=f\volume$ for 
some function $f$ on $M$. By the Cauchy-Schwarz inequality we have 
$$
\int_M\abs{f}\volume\leq\|f\|_{L^2}\sqrt{\text{Vol}(M)},
$$ 
with equality if and only if $f$ is a constant. Hence 
\begin{equation*}
\begin{split}
E(\phi)=&\frac12\|f\|^2_{L^2}\geq
\frac{1}{2\mathrm{Vol}(M)}\big(\int_M \abs{f}\volume\big)^2\\
\geq&\frac{1}{2\mathrm{Vol}(M)}\big(\int_M f\volume\big)^2=
\frac{1}{2\mathrm{Vol}(M)}\big(\int_M\phi^*\omega\big)^2,
\end{split}
\end{equation*}
and the right hand side is a homotopy invariant. This lower bound is 
attained if and only if $f$ is a constant. If $f=0$ then $\phi^*\omega=0$
so $\phi$ is isotropic. If $f\neq 0$ then $\phi$ has no critical points. 
\end{proof}

\begin{example}
Assume that $\phi:M^2 \ra S^2$ is anisotropic and
attains the bound. Then $\phi$ is necessarily a covering map; since $S^2$ 
is simply connected we must have $M=S^2$ and $\phi$ a diffeomorphism. But
then $\int_M\phi^*\omega=\pm\text{Vol}(M)$, so $E(\phi)=\frac{1}{2}
\text{Vol}(M)$. Thus $f\equiv\pm 1$, so $\phi$ is, up to an orientation
reversing isometry, a symplectomorphism of $S^2$. Hence the set of maps
attaining the bound consists of the orbit of $\id:S^2\ra S^2$ under
the group of symplectomorphisms of $S^2$ and the image of this orbit
under reflexion.
\end{example}

\begin{example} \label{t2}
Let $M=N=T^2=\C/(\Z+i\Z)$. Then the bound is attained in
every homotopy class by 
$$
\phi(x,y)=(x,y)L\qquad (L\in\text{Mat}_2(\Z)).
$$
\end{example}

\begin{proposition}\label{bound2} Let $M=M^4$. Then
$$
E(\phi)\geq\frac12\abs{\int_M\phi^*(\omega\wedge\omega)}
$$
with equality if and only if $\phi^*\omega$ is (anti-)self-dual. Note
that the right hand side is a homotopy invariant.
\end{proposition}

\begin{proof}
Since $\dim M=4$, the Hodge $\ast$-operator is an involution on the bundle 
of 2-forms. Thus any 2-form $\alpha$ can be decomposed as 
$$
\alpha=\alpha_++\alpha_-,
$$ 
where $\ast\alpha_+=\alpha_+$ and $\ast\alpha_-=-\alpha_-$. Since 
$\alpha_+$ and $\alpha_-$ are mutually orthogonal we get 
$$
\|\alpha\|^2_{L^2}=\|\alpha_+\|^2_{L^2}+\|\alpha_-\|^2_{L^2}\geq
\abs{\|\alpha_+\|^2_{L^2}-\|\alpha_-\|^2_{L^2}}=
\abs{\int_M\alpha\wedge\alpha},
$$ 
with equality if and only if $\alpha_+=0$ ($\alpha$ is anti-self-dual) or 
$\alpha_-=0$ ($\alpha$ is self-dual). The proposition follows once we
apply this to $\alpha=\phi^*\omega$.
\end{proof}

\begin{remark} The proposition tells us nothing useful if $H^2(M;\R)=0$
or if $\dim N=2$. In the first case $\phi^*\omega$ is necessarily exact,
and hence so is $\phi^*\omega\wedge\phi^*\omega=\phi^*(\omega\wedge\omega)$,
so we just recover the trivial fact that $E(\phi)\geq 0$. In the second
case $\omega\wedge\omega=0$, so again we deduce only $E(\phi)\geq 0$.
\end{remark}

\begin{example} In the case $M=N=T^4$, the bound is attained in each
homotopy class, by linear maps with integer coefficients, as in 
Example \ref{t2}.
\end{example}

Let us now find an explicit formula for the Hessian of a critical point. Note that the metric $h$ on $TN$ induces a metric, also denoted by $h$, on $\phi^{-1}TN$. 
\begin{proposition}\label{hess} Assume that $\phi$ is a critical point of the 
energy functional. Then the Hessian of $\phi$ is given by 
$$
H_{\phi}(X,Y)=\int_M h(X,\L_\phi Y)\volume\qquad(X,Y\in\Gamma(\phi^{-1}TN)),
$$ 
where 
$$
\L_\phi Y=-J\bigg(\nabla^\phi_{Z_\phi}Y+
\d\phi(\musicd\dstar\d\phi^*\iota_Y\omega)\bigg)\ \text{ and }\  Z_\phi=\musicd\dstar\phi^*\omega.
$$
\end{proposition}

\begin{remark} Note that $Z_\phi$ is the vector field on $M$ which must lie
pointwise in $\ker d\phi$ given that $\phi$ is critical, by Corollary 
\ref{cor:ker}.
\end{remark}

\begin{proof}[Proof of Proposition \ref{hess}:] Let $X_s=\partial_t\phi_{s,t}|_{t=0}$. Then, using the 
Homotopy Lemma and a calculation similar to that of the previous section, 
\begin{equation}\label{eq:secondder}
\begin{split}
\frac{\partial^2}{\partial s\partial t}\big|_{s=t=0}\frac12\phi^*_{s,t}
\omega\wedge\ast\phi^*_{s,t}\omega=&\frac{\partial}{\partial s}\big|_{s=0}
\d\phi^*_s\iota_{X_s}\omega\wedge\ast\phi^*_s\omega\\
=&\big(\frac{\partial}{\partial s}\big|_{s=0}\d\phi^*_s
\iota_{X_s}\omega\big)\wedge\ast\phi^*\omega+\d\phi^*\iota_X\omega\wedge
\ast\d\phi^*\iota_Y\omega.
\end{split}
\end{equation}
When integrating over $M$, the second term on the right becomes 
$$
\ip{\d\phi^*\iota_X\omega}{\d\phi^*\iota_Y\omega}_{L^2}=
\ip{\phi^*\iota_X\omega}{\dstar\d\phi^*\iota_Y\omega}_{L^2}.
$$ 
Now note that, using a local orthonormal frame $e_1,\dots,e_m$ for $M$, we 
get  
$$
\ip{\phi^*\iota_X\omega}{\dstar\d\phi^*\iota_Y\omega}=
\sum_{i=1}^m\phi^*\iota_X\omega(e_i)\dstar\d\phi^*\iota_Y\omega(e_i)=
\omega(X,\d\phi(\musicd\dstar\d\phi^*\iota_Y\omega)).
$$ 
Let us now look at the first term on the right hand side of 
\eqref{eq:secondder}. Pointwise we have 
\begin{equation*}
\begin{split}
\ip{\frac{\partial}{\partial s}\big|_{s=0}\d\phi^*_s\iota_{X_s}\omega}{
\dstar\phi^*\omega}=&\sum_{i=1}^m\big(
\frac{\partial}{\partial s}\big|_{s=0}\d\phi^*_s\iota_{X_s}\omega(e_i)\big)
\dstar\phi^*\omega(e_i)\\
=&\big(\frac{\partial}{\partial s}\big|_{s=0}\phi^*_s\iota_{X_s}\omega\big)
(Z_\phi)\\
=&\frac{\partial}{\partial s}\big|_{s=0}\omega(X_s,\d\phi_s(Z_\phi))\\
=&\omega(\nabla^{\phi_s}_{\partial_s}X_s\big|_{s=0},\d\phi(Z_\phi))+
\omega(X,\nabla^{\phi_s}_{\partial_s}\d\phi_s(Z_\phi)\big|_{s=0}).
\end{split}
\end{equation*}
The first term on the right vanishes since we assume $\phi$ to be a 
critical point. The second term becomes 
$$
\omega(X,\nabla^{\phi_s}_{\partial_s}\d\phi_s(Z_\phi)\big|_{s=0})=
\omega(X,\nabla^{\phi_s}_{Z_\phi}\d\phi_s(\partial_s)\big|_{s=0})=
\omega(X,\nabla^\phi_{Z_\phi}Y),
$$ 
where we used the fact that 
$$
\nabla^{\phi_s}_{\partial_s}\d\phi_s(Z_\phi)=\nabla^{\phi_s}_{Z_\phi}
\d\phi_s(\partial_s)+\d\phi_s([\partial_s,Z_\phi])=\nabla^{\phi_s}_{Z_\phi}
\d\phi_s(\partial_s)$$ since $[\partial_s,Z_\phi]=0$. 
\end{proof}

\begin{corollary} Let $\phi:M\to N$ be a critical point of the 
energy functional. Then 
$$
H_\phi(Y,Y)=\int_M \omega(Y,\nabla^\phi_{Z_\phi} Y)\volume\:\:+\:\:
\|\d\phi^*\iota_Y\omega\|^2_{L^2}\qquad(Y\in\Gamma(\phi^{-1}TN)).
$$ 
In particular, $\phi$ is stable if $Z_\phi$ vanishes. 
\end{corollary}
\begin{proof} Take a local orthonormal frame $e_1,\dots,e_m$ for $M$. Then 
\begin{equation*}
\begin{split}
h(JY,\d\phi(\musicd\dstar\d\phi^*\iota_Y\omega))=&\sum_{i=1}^mh(JY,\d\phi(e_i))\dstar\d\phi^*\iota_Y\omega(e_i)\\
=&\sum_{i=1}^m\phi^*\iota_Y\omega(e_i)\dstar\d\phi^*\iota_Y\omega(e_i)\\
=&\ip{\phi^*\iota_Y\omega}{\dstar\d\phi^*\iota_Y\omega}.
\end{split}
\end{equation*}
When integrated, this becomes $\|\d\phi^*\iota_Y\omega\|^2_{L^2}$, and the proof follows. 
\end{proof}

\begin{example} Assume that $\phi:M\to N$ is a critical immersion. 
According to Corollary \ref{cor:ker}, $Z_\phi=0$, so $\phi$ is stable. 
In particular, the identity map of any compact K\"ahler manifold is stable. 
In the case $\dim N=2$ or $4$, we have the stronger 
information that $\id:N\ra N$
globally minimizes $E$ within its homotopy class, by Propositions
\ref{bound1} and \ref{bound2}.

In fact, we can obtain some idea of $\spec\L_\id$, the
eigenvalue spectrum of $\L_\id$, from the following result. Note that
$$
\L_\id=J^{-1}\musicu^{-1}\:\dstar\d\:\musicu J,
$$
so that $Y$ is an eigensection of $\L_\id$ with eigenvalue $\lambda$ if
and only if $\musicu J Y$ is an eigenform of $\dstar\d$ with eigenvalue
$\lambda$.
\begin{proposition} Denote by $\spec\dstar\d$ and $\spec\d\dstar$ the 
eigenvalue spectra of the operators 
$$
\dstar\d:\Gamma(T^*M)\to\Gamma(T^*M)
$$ 
and 
$$
\d\dstar:\Gamma(\wedge^3T^*M)\to\Gamma(\wedge^3T^*M),
$$ 
respectively. Then 
$$
\spec\Delta_2=\spec\dstar\d\cup\spec\d\dstar,
$$
where $\Delta_2$ denotes the Laplacian on $2$-forms on $M$. 
Hence $\spec\L_\id=\spec\dstar\d\subseteq\spec\Delta_2$ in general. 
Furthermore, 
if $\dim M=2$ or $\dim M=4$, then 
\beq\label{spectrum}
\spec\Delta_2=\spec\dstar\d=\spec\L_\id.
\eeq
\end{proposition}
\begin{proof} Since the K\"ahler form on $M$ is harmonic, 
$0\in\spec\Delta_2$. Assume that $\lambda\in\spec\dstar\d$. To show that 
$\lambda\in\spec\Delta_2$, we may thus assume that $\lambda\neq0$. Then 
there is a $1$-form $\alpha\neq0$ on $M$ with 
$$
\dstar\d\alpha=\lambda\alpha.
$$ 
Then we must have $\d\alpha\neq0$ and so 
$$
\d\dstar\d\alpha=\lambda\d\alpha,
$$ 
implying that 
$$
\Delta_2\d\alpha=\lambda\d\alpha.
$$ Hence $\spec\dstar\d\subseteq\spec\Delta_2$. The proof that 
$\spec\d\dstar\subseteq\spec\Delta_2$ is similar. 

Conversely, assume $\lambda\in\spec\Delta_2$. As 
$0\in\spec\dstar\d\cup\spec\d\dstar$, we may assume that $\lambda\neq0$. 
Then there is a $2$-form $\xi\neq0$ with 
\beq\label{jsg1}
\Delta_2\xi=\lambda\xi.
\eeq
By the Hodge decomposition we can write 
$$
\xi=\xi_H+\d\alpha+\delta\beta,
$$ 
where $\xi_H$ is a harmonic $2$-form, $\alpha$ a $1$-form and $\beta$ a 
$3$-form. By the Hodge decomposition of $\alpha$ and $\beta$ we see that we 
may assume that $\alpha$ is coexact and $\beta$ exact. From equation
(\ref{jsg1}) it quickly follows that 
\begin{gather*} \xi_H=0\\
\d\dstar\d\alpha=\lambda\d\alpha\\
\dstar\d\dstar\beta=\lambda\dstar\beta.
\end{gather*} 
The second of these equations implies that the $1$-form 
$\dstar\d\alpha-\lambda\alpha$ is closed; by assumption it is also coexact, 
and so it must vanish. If $\alpha\neq0$ we thus have 
$\lambda\in\spec\dstar\d.$ On the other hand, the third equation implies 
that the $3$-form $\d\dstar\beta-\lambda\beta$ is coclosed; by assumption 
it is also exact, so it must vanish. If $\beta\neq 0$ we thus have 
$\lambda\in\spec\d\dstar.$ Since at least one of $\alpha$ and $\beta$ is 
non-zero, we must have 
$$
\lambda\in\spec\dstar\d\cup\spec\d\dstar.
$$ 

The last statement is obvious when $\dim M=2$ since any $3$-form vanishes. 
When $\dim M=4$,  the action of $\d\dstar$ on $3$-forms is equivalent to 
the action of $\dstar\d$ on $1$-forms under the Hodge isomorphism.   
\end{proof}

It is interesting to compare this with the behaviour
of $\id:(M,g)\ra (M,g)$ in harmonic
map theory, where $\id$ is not stable in general 
(though it is stable if $M$ is
K\"ahler) \cite{lic,smi}. 
The analogous operator to $\L_\id$ is the Jacobi operator
$J_\id$, whose spectral properties depend crucially
 on the Ricci curvature of
$M$. There is a formula similar to (\ref{spectrum}), 
$$
\spec J_\id=\spec\Delta_1-\frac{2s}{\dim M},
$$ 
where $s$ is the scalar curvature of
$M$, but it holds only in
the case that $(M,g)$ is Einstein. More generally there is no
simple relationship between $\spec J_\id$ and $\spec\Delta_p$. 
Analytically, $J_\id$ is elliptic, and so has finite-dimensional kernel.
By contrast, $\ker\L_\id$ is the space of symplectic vector fields
(those $Y$ for which $\iota_Y\omega$ is closed), which has infinite
dimension. Clearly $\dim\ker\L_\phi=\infty$ for all critical maps due to the
invariance of $E(\phi)$ under symplectic diffeomorphisms of $N$.
In fact we will see in the next section an example of a 
critical map $\phi$ (the Hopf map $S^3\ra S^2$)
for which \emph{every}
eigenspace of $\L_\phi$ has infinite dimension.
\end{example}

\begin{example} In Example \ref{ex:warped} we proved that the projection of 
a warped product 
$$
\phi:P\times_f N\to N
$$ 
is critical when $N$ is a K\"ahler manifold. This is not an immersion.
However, we showed that
$Z_\phi=0$, so such a map is always stable nonetheless. The same is
true of the critical projection $(\R^4,g)\ra S^2$ of Example
\ref{funky}.
\end{example}

\section{The Hopf Map}\label{hopf}

In this section we prove that the Hopf map $S^3\ra S^2$
is stable and calculate the spectrum of its Hessian. We then apply this to 
prove a conjecture of Ward regarding the full Faddeev-Hopf model.

We begin by introducing some Lie group and Lie algebra technicalities regarding 
symmetric and Hermitian symmetric spaces. We stringently follow the 
conventions used in \cite{He:78}, to which the reader is referred for 
definitions and fundamental results on symmetric spaces.

Assume that $G$ is a compact, connected, simple Lie group and that $K$ is a 
compact subgroup of $G$ such that $G/K$ is an irreducible Hermitian 
symmetric space of compact type. On the Lie algebra level we have the 
standard orthogonal decomposition 
$$
\g=\k+\p,
$$ 
where $\k$ is the Lie algebra of $K$ and $\p$ an $\Ad_K$-invariant subspace 
with the property that $[\p,\p]\subseteq\k$. It is well known that the 
Hermitian structure on $G/K$ is induced by the adjoint action of an element 
in the centre of $\k$; in accordance with earlier notation, we denote this element by 
$J$.

We provide $G$ with the Riemannian metric induced by the negative of the 
Killing form (or a suitable multiple thereof), and give $G/K$ the metric 
which turns the homogeneous projection 
$$
\phi:G\to G/K,\quad g\mapsto g\cdot o
$$ 
into a Riemannian submersion; here $o$ denotes the identity coset in $G/K$.  
The fibres of $\phi$ are clearly minimal, even totally geodesic; according 
to Example \ref{Example:Riemann}, $\phi$ is a critical point of the 
functional. For simplicity, we denote by $\ep$ the negative of the Killing 
form on $\g$.

The pullback bundle $\phi^{-1}TG/K$ is isomorphic to the trivial bundle 
$G\times\p$ by the map 
$$
G\times\p\ni(g,X)\mapsto\frac{d}{dt}\big|_{t=0}
\phi(g)\exp{tX}\cdot o\in T_{\phi(g)}G/K;
$$ 
the metric on $\phi^{-1}TG/K$ corresponds under this isomorphism to the 
metric on $G\times\p$ induced by the restriction of $\ep$ to $\p$. 
Similarly, we can identify $TG$ with the trivial bundle $G\times\g$ by left 
translation, and this gives the following commutative diagram: 
\begin{equation*}
\xymatrix{ G\times\g\ar[r]^\cong\ar[d] & TG\ar[d]^{d\phi}\\
G\times\p\ar[r]^\cong & \phi^{-1}TG/K}
\end{equation*}
The map on the left, which we thus identify with $\d\phi$, is induced by 
orthogonal projection 
$$
\g=\k+\p\to\p.
$$ 
With this identification in mind, we think of sections of $TG$ as functions 
on $G$ with values in $\g$ and sections of $\phi^{-1}TG/K$ as functions on 
$G$ with values in $\p$: 
$$
\Gamma(TG)\cong C^\infty(G,\g),\quad 
\Gamma(\phi^{-1}TG/K)\cong C^\infty(G,\p).
$$ 
The sections of $TG$ are of course also derivations: for any vector space 
$V$ and smooth function $f:G\to V$, an element $X\in C^\infty(G,\g)$ acts 
on $f$ as 
$$
X(f)(g)=\d f(X)(g)=\frac{d}{dt}\big|_{t=0}f(g\exp(tX(g))\qquad(g\in G).
$$ 
The Levi-Civita connexion on $TG$ corresponds to the connexion 
$$
\nabla_XY=\d Y(X)+\frac12[X,Y]\qquad(X,Y\in C^\infty(G,\g)),
$$ 
and the pullback of the Levi-Civita connexion on $TG/K$ to $\phi^{-1}TG/K$ 
corresponds to the connexion 
$$\nabla^\phi_XY=\d Y(X)+[X,Y]\qquad
(X\in C^\infty(G,\g),\ Y\in C^\infty(G,\p)).
$$

\begin{proposition} Choose an orthonormal basis $\{e_k\}_{k=1}^m$ for $\g$ 
such that $e_1,\dots,e_{m-n}$ is a basis for $\k$ and $e_{m-n+1},\dots,e_m$ 
a basis for $\p$. For the homogeneous projection $\phi$, the second 
variation takes the form 
\begin{equation*}
\begin{split}
\L_\phi Y=-J\bigg(-\frac{\lambda}{2}J(Y)-&J\sum_{k=1}^me_ke_k(Y)+\frac32J
\sum_{a=1}^n[e_a(Y),e_a]\\
+&\sum_{r,s=m-n+1}^m\omega(e_re_s(Y)-\frac12[e_r,e_s](Y),e_r)e_s
\bigg).
\end{split}
\end{equation*}
Here $\lambda$ is the eigenvalue of the Casimir operator associated to the 
adjoint representation of $\g$:
$$
-\sum_{k=1}^m[e_k,[e_k,X]]=\lambda X\qquad(X\in\g).
$$
\end{proposition}
\begin{proof} We begin by calculating $Z_\phi$. By Corollary \ref{cor:ker} 
we know that $\dstar\phi^*\omega(X)=0$ for $X\in\p$, and for $X\in\k$ we have 
\begin{equation*}
\begin{split}
\dstar\phi^*\omega(X)=&\sum_{k=1}^m\big(-e_k(\phi^*\omega(e_k,X))+
\phi^*\omega(\nabla_{e_k}e_k,X)+\phi^*\omega(e_k,\nabla_{e_k}X)\big)\\
=&\sum_{k=1}^m\omega(\d\phi(e_k),\d\phi(\nabla_{e_k}X))\\
=&\sum_{r=m-n+1}^m\omega(\d\phi(e_r),\d\phi(\frac12[e_r,X]))\\
=&\frac12\sum_{r=m-n+1}^m\omega(e_r,[e_r,X])\\
=&\frac12\sum_{r=m-n+1}^m\ip{[Je_r,e_r]}{X}.
\end{split}
\end{equation*}
Thus, $$Z_\phi=\frac12\sum_{r=m-n+1}^m[Je_r,e_r]=\frac12
\sum_{r=m-n+1}^m[e_r,[e_r,J]]=\frac12\sum_{k=1}^m[e_k,[e_k,J]]=
-\frac{\lambda}{2}J,$$ where we have used the fact that $J$ belongs to the centre of 
$\k$. 

Next we look at $\phi^*\iota_Y\omega$. Let $A,B\in\g$. Then 
\begin{equation*}
\begin{split}
\d\phi^*\iota_Y\omega(A,B)=&A(\omega(Y,\d\phi(B))-B(\omega(Y,\d\phi(A))
-\omega(Y,\d\phi([A,B]))\\
=&\omega(\d Y(A),B_\p)-\omega(\d Y(B),A_\p)-\omega(Y,[A,B]_\p).
\end{split}
\end{equation*}
Thus, for $X\in\p$, 
\begin{equation*}
\begin{split}
\dstar\d\phi^*\iota_Y\omega(X)=&\sum_{k=1}^m
\big(-e_k(\d\phi^*\iota_Y\omega(e_k,X))+\d\phi^*\iota_Y
\omega(\nabla_{e_k}e_k,X)\\
&+\d\phi^*\iota\omega(e_k,\nabla_{e_k}X)\big)\\
=&\sum_{k=1}^m\big(-e_k(\omega(\d Y(e_k),X)-\omega(\d Y(X),(e_k)_\p)
-\omega(Y,[e_k,X]_\p)\\
&+\omega(\d Y(e_k),(\nabla_{e_k}X)_\p)-\omega(\d Y(\nabla_{e_k}X),(e_k)_\p)
-\omega(Y,[e_k,\nabla_{e_k}X]_\p)\big)\\
=&\sum_{k=1}^m\big(-\omega(e_ke_k(Y),X)-\frac12
\omega(Y,[e_k,[e_k,X]]_\p)\big)\\
&+\frac32\sum_{a=1}^{m-n}\omega(e_a(Y),[e_a,X])\\
&+\sum_{r=m-n+1}\big(\omega(e_rX(Y),e_r)-\frac12
\omega([e_r,X](Y),e_r)\big)\\
=&\ip{-J\sum_{k=1}^me_ke_k(Y)+\frac{\lambda}{2}JY+\frac32
\sum_{a=1}^{m-n}J[e_a(Y),e_a]}{X}\\
&+\sum_{r=m-n+1}^m\omega((e_rX-\frac12[e_r,X])(Y),e_r).
\end{split}
\end{equation*}
Hence \begin{equation*}
\begin{split}
\d\phi(\musicd\dstar\d\phi^*\iota_Y\omega)=&-J\sum_{k=1}^me_ke_k(Y)+
\frac{\lambda}{2}JY+\frac32\sum_{a=1}^{m-n}J[e_a(Y),e_a]\\
&+\sum_{r,s=m-n+1}^m\omega((e_re_s-\frac12[e_r,e_s])(Y),e_r)e_s.
\end{split}
\end{equation*}
We thus arrive at 
\begin{equation*}
\begin{split}
\L_\phi(Y)=&-J\bigg(-\frac{\lambda}{2}J(Y)-J\sum_{k=1}^me_ke_k(Y)+
\frac32\sum_{a=1}^{m-n}J[e_a(Y),e_a]\\
&+\sum_{r,s=m-n+1}^m\omega((e_re_s-\frac12[e_r,e_s])(Y),e_r)e_s\bigg).
\end{split}
\end{equation*}
\end{proof}
The Hopf map is by definition the map 
$$
\phi:S^3\subset\cn^2\to\cn P^1,\quad \phi(z_1,z_2)=[z_1,z_2].
$$ 
By the identifications 
$$
S^3\cong\SU2,\quad \cn P^1\cong\SU2/\S(\U1\times\U1),
$$ 
we get the alternative definition of the Hopf map as the homogeneous 
projection 
$$
\phi:\SU 2\to\SU 2/\S(\U1\times\U1).
$$ 

Using the previous result and some representation theory, we shall prove 
the following result:
\begin{theorem}\label{the:Hopf} The Hopf map is stable; the Hessian has 
eigenvalues $\dfrac14(n^2+2n)$ and $\dfrac14(n-2k)^2$, 
$k=0,\dots,n$, $n=1,2,\dots$. Each eigenspace is of infinite dimension. 
\end{theorem}
It is well known that, as a harmonic map, the Hopf map is unstable, see \cite{ura}. Returning to the full Faddeev-Hopf model for maps 
$\SU2\ra\S(\U1\times\U1)$,
$$
E(\phi)=\frac{1}{2}\int_{\SU2}(|\d\phi|^2+\alpha|\phi^*\omega|^2)\volume,
$$
we can give precise information on the stability of the Hopf map
for this functional, thus proving a conjecture stated 
by Ward in \cite{war}. The stability properties turn out to be
exactly analogous to those of the identity map in the Skyrme model
on $S^3$ \cite{man}.
\begin{theorem}\label{the:Ward} The Hopf map is an unstable critical
point of the Faddeev-Hopf energy functional if $\alpha<1$ and a stable 
critical point if $\alpha\geq1$. 
\end{theorem} 
The remaining part of this section is devoted to the proof of these two 
results, beginning with Theorem \ref{the:Hopf}.  

The Lie algebra $\su2$ of $\SU2$ has a basis $$\vartheta_1=\frac{i}{2}\begin{pmatrix} 0 & 1 \\ 1 & 0\end{pmatrix},\ \vartheta_2=\frac12\begin{pmatrix}0 & 1 \\ -1 & 0\end{pmatrix},\ \vartheta_3=\frac{i}{2}\begin{pmatrix}1 & 0 \\ 0 & -1\end{pmatrix},$$ and we choose the inner product $\ep$ on $\su2$ which makes this an orthonormal basis; then $\ep$ is a multiple of the Killing form. The isotropy subalgebra $\k=\s(\u1\times\u1)$ is one-dimensional and spanned by $\vartheta_3$, which also acts as $J$ on $\p=\span\{\vartheta_1,\vartheta_2\}$. 

It is easy to see that, for the adjoint representation of $\su2$, the Casimir operator is just multiplication by $2$, i.e., $\lambda=2$. To calculate $\L_\phi$, let $f\in C^\infty(G,\rn)$. Then 
\begin{equation*}
\begin{split}
\L_\phi(f\vartheta_1)=&-J\bigg(-\vartheta_3(f)\vartheta_1-(\vartheta_1^2+\vartheta_2^2+\vartheta_3^2)(f)[\vartheta_3,\vartheta_1]+\frac32\vartheta_3(f)[\vartheta_3,[\vartheta_1,\vartheta_3]]\\
&+(\vartheta_2\vartheta_1-\frac12[\vartheta_2,\vartheta_1])(f)\ip{[\vartheta_3,\vartheta_1]}{\vartheta_2}\vartheta_1+\vartheta_2^2(f)\ip{[\vartheta_3,\vartheta_1]}{\vartheta_2}\vartheta_2\bigg)\\
=&-J\bigg(-\vartheta_3(f)\vartheta_1+(\vartheta_1^2+\vartheta_2^2+\vartheta_3^2)(f)\vartheta_2+\frac32\vartheta_3(f)\vartheta_1\\
&-\vartheta_2\vartheta_1(f)\vartheta_1+\frac12\vartheta_3(f)\vartheta_1-\vartheta_2^2(f)\vartheta_2\bigg)\\
=&-J\bigg((\vartheta_3(f)-\vartheta_2\vartheta_1(f))\vartheta_1+(\vartheta_1^2+\vartheta_3^2)(f)\vartheta_2\bigg)\\
=&(-\vartheta_1^2-\vartheta_3^2)(f)\vartheta_1+(\vartheta_3-\vartheta_2\vartheta_1)(f)\vartheta_2.
\end{split}
\end{equation*}
A similar calculation gives $$\L_\phi(f\vartheta_2)=(-\vartheta_3-\vartheta_1\vartheta_2)(f)\vartheta_1+(-\vartheta_2^2-\vartheta_3^2)(f)\vartheta_2.$$ Hence we can express the differential operator $\L_\phi$ as a matrix using the basis $\{\vartheta_1,\vartheta_2\}$ for $\p$: $$\L_\phi=\begin{pmatrix}-\vartheta_1^2-\vartheta_3^2 & -\vartheta_3-\vartheta_1\vartheta_2 \\ \vartheta_3-\vartheta_2\vartheta_1 & -\vartheta_2^2-\vartheta_3^2\end{pmatrix}.$$

To calculate the spectrum of $\L_\phi$ we recall the Peter-Weyl theorem 
\cite[p17]{petwey}. According to this, $L^2(\SU2,\rn)$ is the orthogonal sum of the finite-dimensional subspaces spanned by matrix elements for the (finite-dimensional) irreducible unitary representations of $\SU2$. Furthermore, these subspaces are invariant under $\L_\phi$. To calculate the spectrum of $\L_\phi$ it is therefore enough to calculate the spectrum of $\L_\phi$ when restricted to these subspaces. Let us therefore momentarily digress for a study of the irreducible representations of $\SU2$. 

As $\SU2$ is the compact real form of $\SLC2$, all irreducible representations of $\SU2$ are obtained by restriction of the irreducible representations of $\SLC2$. For a basis of $\slc2$, let $$X=\begin{pmatrix} 0 & 1 \\ 0 & 0\end{pmatrix},\ Y=\begin{pmatrix}0 & 0 \\ 1 & 0\end{pmatrix},\ H=\begin{pmatrix}1 & 0 \\ 0 & -1\end{pmatrix}.$$ Denote by $V=\cn^2$ the standard representation of $\SLC2$ and by $V^{(n)}=\mathrm{Sym}^n(V)$ the $n^{\rm th}$ symmetric power of $V$, $n=1,2,\dots$. These are precisely the irreducible, finite-dimensional representations of $\SLC2$, and therefore also of $\SU2$. To study the action of $\su2$ on $V^{(n)}$, recall that there is a highest weight vector $v\in V^{(n)}$; let $$v_k=Y^kv\qquad (k=0,1,\dots,n).$$ We adopt the convention that $v_k=0$ for $k<0$ and $k>n$. Then, if $v$ is suitably chosen, $$Hv_k=(n-2k)v_k,\ Xv_k=k(n-k+1)v_{k-1},\ Yv_k=v_{k+1}.$$ Furthermore, $\{v_k\}_{k=0}^n$ is a basis of $V^{(n)}$. A simple calculation gives that 
\begin{equation*}
\begin{split}
\vartheta_1v_k&=\frac{i}{2}\big(k(n-k+1)v_{k-1}+v_{k+1}\big)\\
\vartheta_2v_k&=\frac12\big(k(n-k+1)v_{k-1}-v_{k+1}\big)\\
\vartheta_3v_k&=\frac{i}{2}(n-2k)v_k\\
(-\vartheta_1^2-\vartheta_2^2)v_k&=\frac12(2kn-2k^2+n)v_k\\
(-\vartheta_1^2-\vartheta_2^2-\vartheta_3^2)v_k&=\frac{1}{4}(n^2+2n)v_k.
\end{split}
\end{equation*}

Let us now return to $\L_\phi$ and study its action on $V^{(n)}\otimes\p$, where we think of $\vartheta_k$ as acting by the representation. Since $\vartheta_3=[\vartheta_2,\vartheta_1]$, we can rewrite $\L_\phi$ as 
\begin{equation}\label{eq:L1}
\L_\phi=(-\vartheta_1^2-\vartheta_2^2-\vartheta_3^2)Id+\begin{pmatrix} \vartheta_2^2 & -\vartheta_2\vartheta_1 \\ -\vartheta_1\vartheta_2 & \vartheta_1^2\end{pmatrix}.
\end{equation} 
Let us denote by $A^{(n)}$ the second operator, as acting on $V^{(n)}\otimes\p$. Then, by our earlier calculations, 
\begin{equation}\label{eq:L2}
\L_\phi=\frac14(n^2+2n)Id+A^{(n)}. 
\end{equation}To find the eigenvalues of $\L_\phi$ on $V^{(n)}\otimes\p$, we must find the eigenvalues of $A^{(n)}$. So assume that $\alpha\otimes\vartheta_1+\beta\otimes\vartheta_2\cong(\alpha,\beta)$ is an eigenvector with eigenvalue $\lambda$. Then 
\begin{equation*}
\begin{cases}\vartheta_2^2\alpha-\vartheta_2\vartheta_1\beta&=\lambda\alpha\\ -\vartheta_1\vartheta_2\alpha+\vartheta_1^2\beta&=\lambda\beta\end{cases} \Longrightarrow (\vartheta_1^2+\vartheta_2^2)(\vartheta_2\alpha-\vartheta_1\beta)=\lambda(\vartheta_2\alpha-\vartheta_1\beta).
\end{equation*}
Again, by our earlier calculations, we see that the only possibility is that either $\lambda=0$, or  $$\lambda=\lambda_k=-\frac12(2kn-2k^2+n)\qquad (k=0,\dots,n).$$ Furthermore, when $\alpha=\dfrac{1}{\lambda_k}\vartheta_2v_k$ and $\beta=-\dfrac{1}{\lambda_k}\vartheta_1v_k$, then it is easy to see that $\alpha\otimes\vartheta_1+\beta\otimes\vartheta_2$ is an eigenvector for $A^{(n)}$ with eigenvalue $\lambda_k$. As the linear map $$V^{(n)}\otimes\p\to V^{(n)},\quad (\alpha,\beta)\mapsto\vartheta_2\alpha-\vartheta_1\beta$$ is surjective, the dimension of its kernel equals $\dim V^{(n)}=n+1$. Hence we conclude that \emph{$\L_\phi$, acting on $V^{(n)}\otimes\p$, has the eigenvalue $\dfrac{1}{4}(n^2+2n)$ of multiplicity $n+1$, and the eigenvalues $$\frac14(n^2+2n)-\frac12(2kn-2k^2+n)=\frac14(n-2k)^2\qquad(k=0,\dots,n),$$ each of multiplicity 1}. 

To calculate the spectrum of $\L_\phi$, this time acting as a differential operator on the space spanned by the matrix elements for the representation $V^{(n)}\otimes\p$ of $\SU2$, choose some $\SU2$-invariant inner product $\epn$ on $V^{(n)}$ and define the functions $$\pi_{kl}(g)=(gv_k,v_l)\qquad(g\in\SU2,\ k,l=0,\dots,n).$$ Then, by the invariance of $\epn$,  $$\vartheta(\pi_{kl})(g)=(g\vartheta v_k,v_l)\qquad(g\in\SU2,\ \vartheta\in\su2,\ k,l=0,\dots,n).$$ It follows that, for $k=0,\dots,n$,  $$\frac{1}{\lambda_k}\vartheta_2(\pi_{kl})\otimes\vartheta_1-\frac{1}{\lambda_k}\otimes\vartheta_1(\pi_{kl})\otimes\vartheta_2\qquad(l=0,\dots,n)$$ is an eigenfunction of $\L_\phi$ with eigenvalue $\dfrac14(n-2k)^2$; the corresponding eigenspace is of dimension $n+1$. Furthermore, it is clear that the kernel of $A^{(n)}$, acting on $\span\{\pi_{kl}\}_{k,l=0}^n\otimes\p$, is of dimension $(n+1)^2$; the action of $\L_\phi$ on $\ker A^{(n)}$ is therefore diagonal with eigenvalue $\dfrac14(n^2+2n)$. 

Thus $\L_\phi$ has the eigenvalue spectrum claimed. Further, we
have shown that $\L_\phi$ is non-negative on an $L^2$ orthogonal 
collection of finite-dimensional spaces spanning $L^2$. Hence
$\langle\cdot,\L_\phi\cdot\rangle$ is non-negative on
$L^2$, so $\phi$ is stable. This completes the proof of
 Theorem \ref{the:Hopf}.

\begin{proof}[Proof of Theorem \ref{the:Ward}] Following Urakawa \cite{ura},
we denote by $D$ the Jacobi operator of the Hopf map $\phi$ with respect to 
the Dirichlet energy. The Hessian of $\phi$ with respect to the 
full Faddeev-Hopf functional 
$$
E(\phi)=\frac{1}{2}\int_{\SU2}(|\d\phi|^2+\alpha|\phi^*\omega|^2)\volume
$$ 
is then obviously
$$
\int_{\SU2}h((D+\alpha\L_\phi)X,Y)\volume\qquad(X,Y\in C^\infty(\SU2,\p)),
$$ 
where, as before, 
$$
\p=\span\{\vartheta_1,\vartheta_2\}\subset\su2.
$$ 
We begin by studying the action of $D+\alpha\L_\phi$ on the spaces 
$V^{(n)}\otimes\p$. From \cite{ura} and Theorem 
\ref{the:Hopf} it follows that $D+\alpha\L_\phi$ is positive semi-definite 
for $n\neq1$ and all $\alpha\geq0$. According to \cite[Corollary 8.12]{ura} 
and \eqref{eq:L1} and \eqref{eq:L2} we have on $V^{(1)}\otimes\p$ 
\begin{equation*}
\begin{split}
D+\alpha\L_\phi=&(1+\alpha)(-\vartheta_1^2-\vartheta_2^2-\vartheta_3^2)Id-
2\begin{pmatrix} 0 & \vartheta_3 \\ -\vartheta_3 & 0 \end{pmatrix}+
\alpha A^{(1)}\\
=&\frac{3}{4}(1+\alpha)Id+\begin{pmatrix} \alpha\vartheta_2^2 & -2\vartheta_3-\alpha\vartheta_2\vartheta_1 \\ 2\vartheta_3-\alpha\vartheta_1\vartheta_2 & \alpha\vartheta_1^2 \end{pmatrix}.
\end{split}
\end{equation*}
Let $\{e_1,e_2\}$ be the standard basis for $V^{(1)}\cong\cn^2$. In the basis $\{e_1\otimes\vartheta_1,e_1\otimes\vartheta_2,e_2\otimes\vartheta_1,e_2\otimes\vartheta_2\}$ for $V^{(1)}\otimes\p$ we can express the last term as the matrix 
$$\begin{pmatrix} -\alpha/4 & -i-i\alpha/4 & 0 & 0 \\
i+i\alpha/4 & -\alpha/4 & 0 & 0 \\
0 & 0 & -\alpha/4 & -i-i\alpha/4\\
0 & 0 & i+i\alpha/4 & -\alpha/4\end{pmatrix}.$$
This matrix has eigenvalues $1$ and $-\dfrac{\alpha}{2}-1$, each with multiplicity $2$. On the eigenspace corresponding to the eigenvalue $1$, the operator $D+\alpha\L_\phi$ is obviously positive semi-definite for all $\alpha$, while on the eigenspace corresponding to the eigenvalue $-\dfrac{\alpha}{2}-1$, we have $$D+\alpha\L_\phi=\frac{3}{4}(1+\alpha)Id-(\frac{\alpha}{2}+1)Id=\frac{\alpha-1}{4}Id.$$ We thus conclude that for $\alpha\geq1$, $D+\alpha\L_\phi$ is positive semi-definite on $V^{(n)}\otimes\p$ for all $n$, while for $\alpha<1$, it is negative definite on a $2$-dimensional subspace of $V^{(1)}\otimes\p$. Hence, the operator $D+\alpha\L_\phi$, acting on the space spanned by the matrix elements for the representation $V^{(n)}\otimes\p$ of $\SU2$, is positive semi-definite for all $n$ if $\alpha\geq1$, while it is negative definite on a $4$-dimensional subspace if $\alpha<1$. Theorem \ref{the:Ward} now follows from the Peter-Weyl Theorem. 
\end{proof}

\end{document}